\documentclass[reqno,12pt]{amsart} 
\usepackage{mathrsfs}
\usepackage[margin=1in]{geometry}
\usepackage{amssymb,amsmath,amsthm,amsfonts}
\usepackage{latexsym}
\usepackage{amssymb}
\usepackage{amsmath}
\usepackage{amsfonts}
\usepackage{amsthm}
\usepackage{amscd}
\usepackage{texdraw}
\usepackage{color}
\usepackage[centertags]{amsmath}

\newlength{\defbaselineskip}
\newcommand{\setlinespacing}[1]%
           {\setlength{\baselineskip}{#1 \defbaselineskip}}

\usepackage[pagebackref,colorlinks=true,urlcolor=blue,linkcolor=red,citecolor=red]{hyperref}
\usepackage{epsfig, graphics, graphicx}
\usepackage{subcaption, comment, bm}
\usepackage[percent]{overpic}
\usepackage{float}
\usepackage{amssymb,bm}
\usepackage{amsthm}
\usepackage{color}
\usepackage[]{amsmath}
\usepackage[]{amsfonts}
\usepackage[]{fancyhdr}
\usepackage[]{graphicx,wrapfig}
\usepackage{enumitem}
\usepackage[utf8]{inputenc}
\usepackage[T1]{fontenc}
\usepackage{mathtools}
\usepackage{array}

\definecolor{titlepagecolor}{cmyk}{1,.60,0,.40}

\DeclareFixedFont{\titlefont}{T1}{ppl}{b}{it}{0.5in}

\def\th@plain{%
  \thm@notefont{}
  \itshape 
}
\def\th@definition{%
  \thm@notefont{}
  \normalfont 
}
\makeatother

  \usepackage[pagewise]{lineno}

\DeclareUnicodeCharacter{2212}{-,+}
\usepackage[none]{hyphenat}
\theoremstyle{plain}
\newtheorem{definition}{Definition}[section]
\newtheorem{theorem}[definition]{Theorem}

\newtheorem{lemma}[definition]{Lemma}
\newtheorem{prop}[definition]{Proposition}

\newcommand{\vh}{\textbf{\textit{h}}}
\newcommand{\vx}{\textbf{\textit{x}}}

\newcommand{\vw}{\textbf{\textit{w}}}
\newcommand{\vv}{\textbf{\textit{v}}}

\newcommand{\vf}{\textbf{\textit{f}}}
\newcommand{\vF}{\textbf{\textit{F}}}
\newcommand{\vG}{\textbf{\textit{G}}}
\newcommand{\vg}{\textbf{\textit{g}}}

\newcommand{\vuu}{\textbf{\textit{u}}}

\usepackage{fancyhdr}
\usepackage{lipsum}

\begin{document}
\title[Reconstruction from moment ray transform]{Reconstruction of a vector field and a symmetric $2$-tensor field from the moment ray transforms in $\mathbb{R}^2$}

\author[R. Bhardwaj]{Rahul Bhardwaj}
\address{Department of Mathematics, Indian Institute of Technology, Ropar, Rupnagar-140001, Punjab, India}
\email{bhardwaj161067@gmail.com}

\author[K. B. Solanki]{Karishman B. Solanki}
\email{karishsolanki002@gmail.com, staff.karishman.solanki@iitrpr.ac.in}

		

\begin{abstract} 
We present a technique for recovering a vector field and a symmetric $2$-tensor field, both real-valued and compactly supported in some strictly convex bounded domain with smooth boundary in the Euclidean plane, from the sum of their attenuated moment ray transforms. In addition, we provide a stability estimate for recovering both the vector field and the symmetric $2$-tensor field from the aforementioned ray transform.
\end{abstract}

\maketitle


\section{Introduction} 
There are various types of linear integral operators in integral geometry operating on functions and tensor fields, appearing as mathematical models in wave optics, computerized tomography and emission tomography which can be given by the attenuated ray transforms (ART) or the attenuated moment ray transform (AMRT). These transforms arise in Doppler tomography, Magneto-Acousto-Electrical tomography (MAET), and anisotropic media; see \cite{norton1989tomographic,Jansson1995,kunyansky2023weighted,hendriks2017bragg,sharafutdinov2012integral,natterer2001mathematics} and references therein.
 
The ART can be seen as a weighted ray transform. The ART is an important mathematical tool used in the investigation of several types of inverse problems, and has various applications in the fields of biology and medicine diagnostics, see \cite{budinger1979emission, natterer2001mathematics}, and in the areas of physical optics, specifically wave optics and photometry. The problem of single-photon emission computerized tomography (SPECT) is one of the important problems that is formulated by ART. The AMRTs are generalized for complex-valued absorption coefficient as well as weight functions of polynomial type along with exponential type; refer to these \cite{derevtsov2021generalized, fujiwara2024inversion} for more details. AMRTs of various orders are connected by the means of the application of the linear part of a transport equation. We refer more works related to the ARTS \cite{finch2003attenuated,kuchment2006generalized} and for theoretical results based on the language of the transport equation see \cite{romanov1994conditional, bukhgeim1995inversion} and references therein. The non-attenuated ray transform and the non-attenuated integral moment transform of zeroth order, both are special cases of the ART  and also coincide with the longitudinal ray transform (LRT),  appearing in the linearization of the boundary rigidity problem \cite{sharafutdinov2012integral,stefanov1998rigidity,uhlmann2021travel}. It has many applications in imaging sciences, notably in seismic imaging, ocean imaging, and medical imaging. It is well known that LRT has a large kernel that contains all potential tensor fields that have decay at the boundary of the support; see, for example, \cite{sharafutdinov2012integral}, and thus recovery of the entire tensor field is not possible from LRT data only. The recovery of the solenoidal part of the symmetric $m$-tensor tensor field from the knowledge of its ray transform (0-th moment) has been studied extensively in various settings; see \cite{holman2009weighted,krishnan2019solenoidal,paternain2013tensor,sharafutdinov2012integral,svetov2012reconstruction} and references therein. For the full recovery of the symmetric tensor field, more data is required in the form of combinations of the longitudinal ray transform (LRT), the transverse ray transform (TRT), mixed ray transform (MRT) or from a set of the moment ray transforms, see \cite{Desai2016,kunyansky2023weighted,louis2022inversion,derevtsov2023ray,denisiuk2023iterative,Mishra_2021,krishnan2020momentum,Rohit_Kumar_Mishra_moment,derevtsov2015tomography} and references therein. The moment ray transform is some kind of extension of the ray transform and it was investigated by Sharafutdinov in \cite{sharafutdinov1986problem}. It has been shown that a symmetric $m$-tensor field can be recovered from knowledge of the moment ray transforms up to order $m$, (see \cite[Section 2.12]{sharafutdinov2012integral}).

Recently, the ART  has been investigated by Sadiq, Fujiwara, Tamasan, Derevtsv, and many others; in various settings, we refer to \cite{sadiq2015range, sadiq2016x,omogbhe2024x, omogbhe2025fourier,fujiwara2019fourier, fujiwara2024inversion} and references therein. In \cite{fujiwara2024inversion}, Fujiwara et al. studied the moment ray transform and provided an approach for full recovery of the symmetric m-tensor field from the knowledge of its attenuated moment ray transforms of order $0$ to $m$, and in \cite{derevtsov2021generalized}, the authors investigated various properties of ART and their integral angular moments. A technique to recover the solenoidal parts of a vector field and a symmetric 2-tensor field from the sum of their ART has been obtained by Omogbhe in \cite{omogbhe2025fourier}. This type of problem comes after linearization of a problem that is related to travel time tomography, described via Mane's action potential of the energy level $1/2$ for magnetic flow \cite{dairbekov2007boundary,stefanov2024lorentzian}.

Motivated by the above works, in this article, we study the attenuated moment ray transform defined by the sum of the $k$-th attenuated moments of the vector field and the 2-tensor field, for $0\leq k \leq 2$, and consider the inverse problem of recovering a vector field and a symmetric 2-tensor field, whose support is contained in the strictly convex bounded domain having smooth boundary in the plane. The idea is to convert the problem into a boundary value problem (BVP) for a system of transport equations, and then solve it using an extension of the Bukhgeim A-analytic theory. We refer the reader to \cite{bukhgeim1995inversion, omogbhe2025fourier,fujiwara2019fourier, fujiwara2024inversion} for works where such technique has been used.

This article is organized as follows. In Section \ref{sec;Preliminaries and statement of the main result}, the notation, definitions and statement of the main theorem are stated. In Section \ref{sec:TM and BA}, we introduce a BVP required to prove our theorem and provide details about the Bukhgeim A-analytic theory. Section \ref{sec;Proof of the main result} is devoted to the proofs, and this section splits into three subsections, \ref{sec;non-attenuated case} and \ref{sec; attenuated case} for the non-attenuated case and  for the attenuated case, respectively, and \ref{sec,proof of Prop} for the proof of a proposition which relates recovery problem with the BVP.


\section{Preliminaries and statement of the main result}\label{sec;Preliminaries and statement of the main result}
In this section, we give some definitions and notations which will be used throughout the article. Let $\Omega$ be a strictly convex bounded domain in $\mathbb{R}^2$ with smooth boundary $\Gamma$, and $\mathbb{S}^1:=\{u_\varphi=(\cos(\varphi),\sin(\varphi)):0\leq\varphi<2\pi\}$ be the unit circle in $\mathbb{R}^2$. 
                        Let $S^{1}(\Omega)$ and $S^{2}(\Omega)$ be the vector spaces of real-valued vector fields and real-valued symmetric $2$-tensor fields supported in $\overline{\Omega}=\Omega\cup\Gamma$, respectively. For $m\in\mathbb{N}$ and $0< \alpha \leq 1$,   $C^{m, \alpha }(\Omega)$ denotes the H\"older space, the space of real-valued functions defined on $\Omega$ that are $m$ times continuously differentiable and whose $m$-th derivatives satisfy the H\"older condition with index $\alpha$. 

For $s\geq0$, $H^s(\mathbb{R}^2)=\{f\in L^2(\mathbb{R}^2): (1+|\xi|^2)^\frac{s}{2} \widehat{f}(\xi)\in L^2(\mathbb{R}^2) \}$ denotes the standard Sobolev space with the norm \begin{align*} 
    \|f\|_{H^s(\mathbb{R}^2)}= \left( \int_{\mathbb{R}^2} |\widehat{f}(\xi)|^2 (1+|\xi|^2)^s \,d\xi \right)^\frac{1}{2},
\end{align*}
where $\widehat{f}$ denotes the Fourier transform of $f$. The space $H_0^{s}(\Omega)$ denotes the Sobolev space which is closure of $C_c^\infty(\Omega)$ in $H^s(\Omega)=\{f\in L^2(\Omega): f=g_{|\Omega} \ \text{for some} \ g\in H^s(\mathbb{R}^2) \}$ with the quotient norm \begin{align*}
    \|f\|_{H^s(\Omega)}=\inf \{ \|g\|_{H^s(\mathbb{R}^2)} : g\in H^s(\mathbb{R}^2), \ g_{|\Omega}=f \}.
\end{align*}

For $m\in\{1,2\}$, let $H_0^{s}(\Omega; S^{m})$ be the space of $S^{m}(\Omega)$ valued functions on $\Omega$ with each component in the Sobolev space $H_0^{s}(\Omega)$, i.e., a vector field ${\vf} \in H_0^{s}(\Omega; S^{1})$ and a symmetric 2-tensor field ${\vF} \in H_0^{s}(\Omega; S^{2})$ are of the form
\begin{align}
{\vf}(\vx) = ({f}_1 (\vx), {f}_2 (\vx)) \quad \text{and} \quad 
{\vF}(\vx) = \begin{bmatrix}
        {F}_{11}(\vx) & {F}_{12}(\vx)\\
        {F}_{12}(\vx) & {F}_{22}(\vx)
    \end{bmatrix} \quad (\vx\in{\Omega}),
\end{align}
where $f_i \in H_0^{s}(\Omega)$ and ${F}_{ij} \in H_0^{s}(\Omega)$ for $1\leq i,j \leq 2$.
We denote by $\langle \cdot, \cdot \rangle$ the standard dot product in $S^2(\Omega)$ given by $\displaystyle \langle \vF,\vG \rangle = \sum\limits_{i,j = 1}^{2}F_{ij}G_{ij}$, for $\vF,\vG \in S^{2}(\Omega)$, and the corresponding norm will be denoted by $|\cdot|$.

\begin{definition}\label{eq:tensor product}
    The \emph{tensor product} of two vectors $\vuu = (u_1, u_2)$ and $\vv = (v_1,v_2)$ in $\mathbb{R}^2$ is a $2$-tensor field, denoted by $\vuu \otimes \vv$ and defined component-wise as  
    \begin{align}\label{tensor product}
        (\vuu\otimes\vv)_{ij} := u_iv_j . 
    \end{align}
     The \emph{symmetrized tensor product} $\vuu \odot \vv$ is defined as 
    \begin{align}\label{symmetrized tensor product}
        \vuu \odot \vv = \frac{1}{2} (\vuu\otimes\vv + \vv\otimes\vuu) .
    \end{align}
\end{definition}

In particular $\vuu^{2}$, the symmetrized tensor product of a vector $\vuu$ with itself, is given by 
    \begin{align}\label{square}
        \vuu^{2} = \vuu \odot \vuu = \vuu \otimes \vuu = \begin{bmatrix}
        \vuu_1^2 & \vuu_1 \vuu_2\\
        \vuu_1 \vuu_2 &  \vuu_2^2
    \end{bmatrix}.
    \end{align}

For ${\vF} \in H_0^{s}(\Omega; S^{2})$ and $\vuu = (u_1, u_2) \in \mathbb{R}^2$, using \eqref{square}, we have the following identity
\begin{align}
    \left\langle {\vF},\vuu^2 \right\rangle = u_{1}^{2}{F}_{11} + 2 u_1u_2{F}_{12} + u_{2}^2{F}_{22}.
\end{align}

Now, we provide the definition of the  $\mathfrak{a}$-attenuated moment ray transform.
\begin{definition}\label{def:moment ray transform}
    Let ${\vf} \in H_0^{s}(\Omega; S^{1})$, ${\vF} \in H_0^{s}(\Omega; S^{2})$, and let $\mathfrak{a} \in C^{m,\alpha}(\Omega)$. For $k\in\{0,1,2\}$, the \emph{$\mathfrak{a}$-attenuated moment ray transforms of order $k$}, denoted by $\mathcal{M}^{(k)}_\mathfrak{a}$, are defined by 
    \begin{align} 
      \mathcal{M}_{\mathfrak{a}}^{(k)}(\mathcal{F})(\vx, \vuu_{\varphi})  := \int\limits_{-\infty}^{\infty} s^{k} \mathcal{F}(\vx +s\vuu_\varphi) \exp \left\{ -\int\limits_{s}^{\infty}\mathfrak{a}(\vx + t \vuu_{\varphi})\,dt\right\} \,ds,
    \end{align}
    where $\mathcal{F}(\vx + s\vuu_\varphi) := {\vf}(\vx +s\vuu_\varphi)\cdot\vuu_{\varphi}+\left\langle {\vF} (\vx + s\vuu_\varphi), \vuu_\varphi^{2} \right\rangle$.
\end{definition}
When $\mathfrak{a} = 0$, the above definition will give \emph{non-attenuated moment ray transforms} denoted by $\mathcal{M}^{(k)}(\mathcal{F})$ $(0\leq k \leq 2)$. In the attenuated case, we assume that the attenuation coefficient $\mathfrak{a} > 0$ in $\overline{\Omega}$.

As done in \cite{denisiuk2023iterative} for non-attenuated case and in \cite{fujiwara2024inversion} for attenuated case, to avoid excessiveness of parametrization $\vx+t\vuu_\varphi$ ($t\in\mathbb{R}$) of line, we shall use a different parametrization of line given by $\Pi_{\varphi}(\vx)+t\vuu_\varphi$ with additional condition $\Pi_{\varphi}(\vx)\cdot\vuu_\varphi=0$. Thus, we give an equivalent definition as
    \begin{align}\label{eq:moment ray transform}
      \mathcal{M}_{\mathfrak{a}}^{(k)}(\mathcal{F})(\vx,\vuu_\varphi)  := \int\limits_{-\infty}^{\infty} s^{k} \mathcal{F}(\Pi_{\varphi}(\vx) +s\vuu_\varphi) \exp \left\{ -\int\limits_{s}^{\infty}\mathfrak{a}(\Pi_{\varphi}(\vx) + t \vuu_{\varphi})\,dt\right\} \,ds,
    \end{align}
where $\Pi_{\varphi}(\vx) := \vx - (\vx\cdot\vuu_{\varphi})\vuu_{\varphi} $ is the projection map onto $\vuu_{\varphi}^{\perp}=(-\sin(\phi),\cos(\phi))$.

Since one of the key approaches in our proof relies on the Fourier series expansion of the functions on $\Omega\times\mathbb{S}^1$, therefore we give brief details about the Fourier series of functions. For fixed $\vx\in\Omega$, let $\displaystyle \vv(\vx, \vuu_{\varphi}) = \sum\limits_{n=-\infty}^{\infty} v_n(\vx)\mathrm{e}^{\mathrm{i}n\varphi}$ be the Fourier series of $\vv(\vx,\cdot)$, where $$v_{-n}(\vx) = \frac{1}{2\pi}\int\limits_{0}^{2\pi} \vv(\vx, \vuu_{\varphi})\mathrm{e}^{in\varphi}\,d\varphi \quad (n\in\mathbb{Z})$$ are the Fourier coefficients of $\vv(\vx,\cdot)$. So, for each $n\in\mathbb{Z}$, $v_{-n}$ becomes a function on $\Omega$. Moreover, if $\vv$ is a real-valued function, then its Fourier coefficients are related via the complex conjugates, i.e., $\overline{v_n} = v_{-n}$, and hence it is enough to work with the sequences of non-positive indices. Define a sequence valued map on $\Omega$ as 
\begin{align} \label{def,seqvalmap}
    \Omega \ni \vx \to \vv(\vx):= (v_0(\vx),v_{-1}(\vx),v_{-2}(\vx),\dots),
\end{align}
and consider the $H^{q}(\Omega)$ valued weighted $\ell^2$-space, for $1\leq p<\infty$ and $q>0$,
\begin{align} \label{def,weiellspa}
    \ell^{2,p}(\mathbb{N}\cup\{0\}, H^q(\Omega)) = \left\{\vv=(v_0,v_{-1},v_{-2},\dots) : \|\vv\|^2_{p,q} = \sum_{k=0}^\infty (1+k)^{2p} \|v_{-k}\|^2_{H^q(\Omega)} \right\}, 
\end{align} 
where the indices $p$ and $q$ are for the smoothness in the angular and the spatial variable, respectively. 

Throughout the article, the notation $\left\lVert \vv\right\rVert \lesssim \left\lVert \vw\right\rVert$ and $\left\lVert \vv\right\rVert \approx \left\lVert \vw\right\rVert$ means $\left\lVert \vv\right\rVert \leq c \left\lVert \vw\right\rVert $ and $c^{-1}\left\lVert \vv\right\rVert\leq \left\lVert \vw\right\rVert\leq c\left\lVert \vv\right\rVert $ for some constant $c>0$, independent of both $v$ and $w$.

In \cite{omogbhe2025fourier}, Omogbhe showed that solenoidal parts of a vector field and a symmetric two tensor field can be recovered from the AMRT \eqref{eq:moment ray transform} of order 0. Since Omogbhe only considered the AMRT of order $0$, and we know that this transform of a vector field and a symmetric 2-tensor field has a non-trivial kernel, therefore it is not possible to recover them fully using only the ray transform data. So, we need more data in order to fully recover them. Here we provide a technique for the full recovery of the vector field and the symmetric 2-tensor field from a set of AMRTs of orders $0$ to $2$ as in \eqref{eq:moment ray transform}. In addition we provide the stability estimate for the aforementioned AMRTs and it is our main theorem which we state now.

\begin{theorem}[Main Theorem] \label{maintheorem}
Let $\Omega$ be a strictly convex bounded domain in $\mathbb{R}^2$ with the smooth boundary $\Gamma$. If ${\vf} \in H_0^{s}(\Omega; S^{1})$, ${\vF} \in H_0^{s}(\Omega; S^{2})$ and $\mathfrak{a} \in C^{1,\alpha}(\overline{\Omega})$, for $s\geq0$ and $1/2<\alpha\leq1$, then $\vf$ and $\vF$ can be recovered uniquely from the $\mathfrak{a}$-attenuated moment ray transforms $\mathcal{M}_{\mathfrak{a}}(\mathcal{F}) := (\mathcal{M}_{\mathfrak{a}}^{(0)}(\mathcal{F}), \mathcal{M}_{\mathfrak{a}}^{(1)}(\mathcal{F}), \mathcal{M}_{\mathfrak{a}}^{(2)}(\mathcal{F}))$. \\ Moreover, if ${\vf} \in H_0^\frac{5}{2}(\Omega; S^{1})$ and ${\vF} \in H_0^\frac{7}{2}(\Omega; S^{2})$, then $\mathcal{M}_{\mathfrak{a}}^{(k)}(\mathcal{F}) \in H^\frac{7}{2}(\mathbb{S}^1,H^\frac{5}{2}(\Gamma))$, for $0\leq k\leq 2$, and the following stability estimate holds
\begin{align} \label{stability estimate}
    \left\lVert \vf\right\rVert_{L^2(\Omega)} + \left\lVert\vF\right\rVert_{L^2(\Omega)} \lesssim \sum_{k=0}^2 \left\lVert\mathcal{M}_{\mathfrak{a}}^{(k)}(\mathcal{F})\right\rVert_{\frac{7}{2},k+\frac{1}{2}}
\end{align} 
\end{theorem}

To prove  the above theorem, we first  introduce a BVP for  the transport equation having  the values of transforms on a part of the boundary and then we recover the vector field and the 2-tensor field, with the help of  Bukhgeim A-analytic theory as used in \cite{omogbhe2025fourier, fujiwara2024inversion}.

\section{Transport model and Bukhgeim A-analytic theory} \label{sec:TM and BA}
In this section, we introduce the BVP for the system of transport equations and give some details about Bukhgeim A-analytic theory as per our requirement.

For $(\vx,\vuu_{\varphi}) \in \overline{\Omega}\times \mathbb{S}^1$, $\tau(\vx,\vuu_{\varphi})=\tau_+(\vx,\vuu_{\varphi}) + \tau_-(\vx,\vuu_{\varphi})$ denotes the length of the chord that passes through the point $\vx$ in the direction of $\vuu_\varphi$, where $(+)$ is for the outgoing direction and $(-)$ is for the incoming direction. Correspondingly, define the sub-bundles restricted to the boundary as 
\begin{align}
    \Gamma_{\pm} &:= \left\{(\vx,\vuu_{\varphi})\in \Gamma\times\mathbb{S}^1| \pm \vuu_{\varphi}\cdot\nu(\vx) > 0\right\},
\end{align}
where $\nu(\vx)$ is the outward unit normal at the point $\vx$ on the boundary $\Gamma$.

For $0\leq k \leq 2$, the $\mathfrak{a}$-attenuated moment ray transform of $\mathcal{F}$, $\mathcal{M}_{\mathfrak{a}}^{(k)}(\mathcal{F})$, is realized as a function on $\Gamma_+$ by
\begin{align}
     \mathcal{M}_{\mathfrak{a}}^{(k)}(\mathcal{F})(\vx,\varphi)  := \int\limits_{-\tau(\vx,\vuu_{\varphi})}^{0}t^{k} \mathcal{F}(\Pi_{\varphi}(\vx) +s\vuu_\varphi) \exp \left\{ -\int\limits_{s}^{0}\mathfrak{a}(\Pi_{\varphi}(\vx) + t \vuu_{\varphi})\,dt\right\} \,ds.
\end{align}
The reason that the limits of integration are changed in the above integration is that the vector field and the symmetric tensor field in consideration are supported in $\overline{\Omega}$, and the line segment inside $\overline{\Omega}$ is given by these limits.

The next step is to transfer the given data to the BVP for the system of transport equations, which will be proved in the following proposition. Its proof is provided in section \ref{sec,proof of Prop}.

\begin{prop}\label{proposition; transport model}
    Let ${\vf} \in H_0^{s}(\Omega; S^1))$, ${\vF} \in H_0^{s}(\Omega; S^2)$, and let $\mathfrak{a} \in C^{1,\alpha}(\overline{\Omega})$, for some $s\geq0$ and $1/2<\alpha\leq1$. Then the following system of BVP on $\overline{\Omega}\times\mathbb{S}^1$
    \begin{align}\label{eq:system of BVP}
    \begin{rcases}
        \vuu_{\varphi}\cdot\nabla \vv^{0}(\vx,\vuu_{\varphi}) + \mathfrak{a}(\vx) \vv^{0}(\vx,\vuu_{\varphi}) = \mathcal{F}(\vx,\vuu_{\varphi}) :=  {\vf}(\vx)\cdot\vuu_{\varphi}+\left\langle {\vF} (\vx), \vuu_\varphi^{2} \right\rangle\\
        \vuu_{\varphi}\cdot\nabla \vv^{1}(\vx,\vuu_{\varphi}) + \mathfrak{a}(\vx) \vv^{1}(\vx,\vuu_{\varphi}) = \vv^{0}(\vx,\vuu_{\varphi})\\
        \vuu_{\varphi}\cdot\nabla \vv^{2}(\vx,\vuu_{\varphi}) + \mathfrak{a}(\vx) \vv^{2}(\vx,\vuu_{\varphi}) = \vv^{1}(\vx,\vuu_{\varphi})
    \end{rcases}
\end{align}
subject to
\begin{align}
    \vv^{k}|_{\Gamma_-} = 0 \quad (0\leq k \leq 2),
\end{align}
has a unique solution $\vv^k \in H^s(\Omega\times\mathbb{S}^1)$. In particular, if $s\geq 1$, then $\vv^k|_{\Gamma\times\mathbb{S}^1} \in H^s(\mathbb{S}^1;H^{s-\frac{1}{2}}(\Gamma))$. Moreover, $(\vv^{0}|_{\Gamma_{+}}, \vv^{1}|_{\Gamma_{+}}, \vv^{2}|_{\Gamma_{+}})$ and $(\mathcal{M}_{\mathfrak{a}}^{(0)}\mathcal{F}, \mathcal{M}_{\mathfrak{a}}^{(1)}\mathcal{F}, \mathcal{M}_{\mathfrak{a}}^{(2)}\mathcal{F})$ are in a one-to-one correspondence by
\begin{align}\label{system of transport eq.}
    \begin{rcases}
        \vv^{0}|_{\Gamma_{+}}(\vx,\vuu_{\varphi}) = \mathcal{M}_{\mathfrak{a}}^{(0)}\mathcal{F}(\vx,\vuu_{\varphi}),\\
         \vv^{1}|_{\Gamma_{+}}(\vx,\vuu_{\varphi}) = (\vx\cdot\vuu_{\varphi})\vv^{0}|_{\Gamma_{+}}(\vx,\vuu_{\varphi}) - \mathcal{M}_{\mathfrak{a}}^{(1)}\mathcal{F}(\vx,\vuu_{\varphi}),\\
          \vv^{2}|_{\Gamma_{+}}(\vx,\vuu_{\varphi}) = \sum\limits_{n=1}^{2}(-1)^{n-1}\frac{(\vx\cdot\vuu_{\varphi})^n}{n!} \vv^{2-n}|_{\Gamma_{+}}(\vx,\vuu_{\varphi}) + \frac{1}{2}\mathcal{M}_{\mathfrak{a}}^{(2)}\mathcal{F}(\vx,\vuu_{\varphi}).
    \end{rcases}
\end{align}
\end{prop}

Here, we considered the BVP to recover the solution of the system of equations \eqref{eq:system of BVP} together with the unknown vector field ${\vf}$ and the 2-tensor field ${\vF}$ from the knowledge of $\vv^k|_{\Gamma\times\mathbb{S}^1}$, for all $k= 0,1,2$. We refer to $\vv^k(\vx,\vuu_{\varphi})$ in \eqref{eq:system of BVP} as the $k-$level flux for $k= 0,1,2$. To solve this BVP, the A-analytic theory of Bukhgeim \cite{bukhgeim1995inversion} will be used, which is based on the Cauchy problem for a Beltrami-like equation associated with A-analytic maps in the sense of Bukhgeim. In this case, first, we expand the solution of the above system of transport equations into the Fourier series and study its Fourier coefficient by some interpretation of harmonic analysis. For more details, we refer the reader to \cite{arbuzov1998two} and \cite{bukhgeim1995inversion} for the Bukhgeim A-analytic theory.

The following are some notions and important results from harmonic analysis, singular integrals, and Bukhgeim $A$-analytic theory which will be used here.

Let $z=x_1 + \mathrm{i} x_2$. The Cauchy-Riemann operators are given by
\begin{align}
    \overline{\partial} = \frac{\partial_{x_1} +\mathrm{i}\partial_{x_2}}{2}, \quad \partial = \frac{\partial_{x_1} -\mathrm{i}\partial_{x_2}}{2},
\end{align}
and the advection operator is
\begin{align}\label{eq;advection}
  \vuu_{\varphi}\cdot\nabla = \mathrm{e}^{-\mathrm{i}\varphi}\overline{\partial} +  \mathrm{e}^{\mathrm{i}\varphi}\partial.
  \end{align}

Denote by $\ell^1$ and $\ell^\infty$ the spaces of summable and bounded sequences respectively. Let $\vw \in C(\overline{\Omega};\ell^{\infty}) \cap C^1({\Omega};\ell^{\infty})$. Consider $\Omega$ as a subset of $\mathbb{C}$. The sequence-valued map  $\Omega \owns z \mapsto \vw(z):= (w_{0}(z), w_{-1}(z), w_{-2}(z),\dots)$ is \emph{$\mathcal{L}^2$-analytic} (in the sense of Bukhgeim) if it satisfies the homogeneous Beltrami-like equation
    \begin{align}\label{eq;homo Beltrami like eq}
        \overline{\partial}\vw(z) + \mathcal{L}^2\partial \vw(z) = 0 \quad (z\in \Omega),    \end{align}
        where $\mathcal{L}$ is the left shift operator defined by $$\mathcal{L}(w_{0}(z), w_{-1}(z), w_{-2}(z),\dots ) = ( w_{-1}(z), w_{-2}(z),\dots)$$ and $\mathcal{L}^2 = \mathcal{L}\circ \mathcal{L}$. As shown in the original paper of Bukhgeim \cite{bukhgeim1995inversion}, the solution of the above equation \eqref{eq;homo Beltrami like eq} is given by a Cauchy-like integral formula
        \begin{align}\label{eq; solution of homo Beltrami like eq}
            \vw(z) = \mathfrak{B}[\vw|_\Gamma](z) \quad (z\in \Omega),
        \end{align}
        where $\mathfrak{B}$ is the Bukhgeim-Cauchy operator which operates on $\vw|_{\Gamma}$ and is described component wise, for $n\geq 0$ and $z\in \Omega$, as follows (refer to \cite{finch2003attenuated}):
    \begin{align}\label{Bukhgeim-Beltrami operator}
        (\mathfrak{B}\vw)_{-n}(z) := \frac{1}{2\pi \mathrm{i}}\int\limits_{\Gamma} \frac{w_{-n}(\zeta)}{\zeta - z}\,d\zeta  +  \frac{1}{2\pi \mathrm{i}}\int\limits_{\Gamma} \left( \sum\limits_{j=1}^{\infty}w_{-n-2j} \left( \frac{\overline{\zeta} - \overline{z}}{\zeta - z}\right)^{j} \right) \left( \frac{d\zeta}{\zeta - z} - \frac{d\overline{\zeta}}{\overline{\zeta} - \overline{z}} \right). 
    \end{align}

The inhomogeneous Bukhgeim-Beltrami equation is of the form
\begin{align}\label{eq;nonhom_}
    \overline{\partial}\vw(z) + \mathcal{L}^2\partial \vw(z) = \vh(z) \quad (z\in \Omega),
\end{align}
and its solution is obtained with the help of a Pompeiu-like operator $\mathcal{T}$ which is described component wise for $n\geq 0$ as follows: 
\begin{align}\label{Pompeiu-like operator}
    (\mathcal{T}\vh)_{-n}(z) := -\frac{1}{\pi}\sum\limits_{j=0}^{\infty}\int\limits_{\Omega}h_{-n-2j}(\zeta)\frac{1}{\zeta - z}\left( 
 \frac{\overline{\zeta} - \overline{z}}{\zeta - z}\right)^{j}\,d\xi\,d\eta \quad (\zeta = \xi + \mathrm{i}\eta, z\in \Omega).
\end{align}

Lastly, we state a proposition that is required here and its proof can be found in \cite{fujiwara2024inversion}.

\begin{prop} \cite[Proposition B.1]{fujiwara2024inversion} \label{prop_non-homogeneous solution}
    Let $\Omega$ be a bounded convex domain with a $C^1$ boundary, and let $h\in C(\overline{\Omega};\ell^1)$. If $ \vw\in C(\overline{\Omega};\ell^1)\cap C^1({\Omega};\ell^1)$ is a solution of the inhomogeneous Bukhgeim-Beltrami equation \eqref{eq;nonhom_}, then 
    \begin{align}
        \vw(z) = \mathfrak{B}[\vw|_{\Gamma}](z) + (\mathcal{T}\vh)(z) \quad (z\in \Omega),
    \end{align}
     where the operators $\mathfrak{B} $ and $\mathcal{T}$ are given by equations \eqref{Bukhgeim-Beltrami operator} and \eqref{Pompeiu-like operator}, respectively. 
\end{prop}

\section{Proofs} \label{sec;Proof of the main result}
First, we give some necessary details and then prove the main theorem in Subsections \ref{sec;non-attenuated case} and \ref{sec; attenuated case} for the non-attenuated and the attenuated cases, respectively. The proof of Proposition \ref{proposition; transport model} is provided in Subsection \ref{sec,proof of Prop}.

If $\vx=(x_1,x_2)\in\mathbb{R}^2$, consider $z=x_1 + \mathrm{i} x_2 \in \mathbb{C}$ as earlier. We shall use this relation between $z$ and $\vx$ interchangeably for convenience. For $\vuu_{\varphi} = (\cos{\varphi}, \sin{\varphi})$, the expression $\displaystyle \cos \varphi = \frac{e^{\mathrm{i}\varphi} + e^{-\mathrm{i}\varphi}}{2}$ and  $\displaystyle \sin \varphi = \frac{e^{\mathrm{i}\varphi} - e^{-\mathrm{i}\varphi}}{2i}$ gives
\begin{align}
    \left\langle {\vF} (\vx), \vuu_\varphi^{2} \right\rangle &= \mathscr{F}_{0}(z) + \overline{\mathscr{F}_{2}(z)}e^{\mathrm{i}2\varphi} + \mathscr{F}_{2}(z)e^{-\mathrm{i}2\varphi} \label{eq;relation btw components of the function}\\ 
    \text{and} \quad {\vf}(\vx)\cdot\vuu_{\varphi} &=  \overline{\mathscr{F}_{1}(z)}e^{\mathrm{i}\varphi} + \mathscr{F}_{1}(z)e^{-\mathrm{i}\varphi}, \label{eq:relation btw components of the function}
\end{align}
where 
\begin{align}\label{equation;relation btw components of the function}
  \mathscr{F}_{0}(z) = \frac{{F}_{11}(\vx) + {F}_{22}(\vx)}{2}, \ \mathscr{F}_{1}(z) = \frac{{f}_{1}(\vx) + \mathrm{i}{f}_{2}(\vx)}{2}
\  \text{and} \ \mathscr{F}_{2}(z) = \frac{{F}_{11}(\vx) - {F}_{22}(\vx)}{4} + \mathrm{i} \frac{{F}_{12}(\vx)}{2}.
\end{align}
Using the above equations \eqref{eq;relation btw components of the function}, \eqref{eq:relation btw components of the function} 
and \eqref{equation;relation btw components of the function} in the system of the transport equation \eqref{eq:system of BVP}, we get
\begin{align}\label{eq;transport model, after putting the value of F}
\begin{rcases}
        \vuu_{\varphi}\cdot\nabla \vv^{0}(z,\vuu_{\varphi}) + \mathfrak{a}(z) \vv^{0}(z,\vuu_{\varphi}) = \mathscr{F}_{0}(z) + \sum_{j=1}^2 \left(\overline{\mathscr{F}_{j}(z)} e^{\mathrm{i}j\varphi} + \mathscr{F}_{j}(z)e^{-\mathrm{i}j\varphi}\right) \\
        \vuu_{\varphi}\cdot\nabla \vv^{1}(z,\vuu_{\varphi}) + \mathfrak{a}(z) \vv^{1}(x) = \vv^{0}(z,\vuu_{\varphi})\\
        \vuu_{\varphi}\cdot\nabla \vv^{2}(z,\vuu_{\varphi}) + \mathfrak{a}(z) \vv^{2}(z,\vuu_{\varphi}) = \vv^{1}(z,\vuu_{\varphi})
\end{rcases}
\end{align}
subject to
\begin{align}\label{eq;traces}
        \vg^k(z,\vuu_{\varphi}) := \begin{cases}
            \vv^k|_{\Gamma_+}(z,\vuu_{\varphi}), & (z,\vuu_{\varphi})\in \Gamma_{+}\\
            0,  & (z,\vuu_{\varphi})\in \Gamma_{-}
        \end{cases} 
         \quad (0\leq k \leq 2).
\end{align}
In the above problem, the solution $\vv^{k}$ ( $0\leq k \leq 2$) of BVP \eqref{eq;transport model, after putting the value of F} and \eqref{eq;traces} is unknown in $\Omega$ as the symmetric 2-tensor field $\vF$ and the vector field $\vf$ are unknown but their traces $g^k$ are known on the boundary $\Gamma \times \mathbb{S}^1$ by \eqref{system of transport eq.} and \eqref{eq;traces}. With this foundation, we now provide the proof for the non-attenuated case.

\subsection{The non-attenuated case:} \label{sec;non-attenuated case}
In the non-attenuated case, the solution of the BVP \eqref{eq;transport model, after putting the value of F}--\eqref{eq;traces} is denoted by $\vw$. Using the Fourier series expansion $ \vw(z, \vuu_{\varphi}) = \sum\limits_{n=-\infty}^{\infty} w_n(z)\mathrm{e}^{\mathrm{i}n\varphi}$ and the advection equation \eqref{eq;advection} into the above BVP \eqref{eq;transport model, after putting the value of F}--\eqref{eq;traces}, we have the following relations for the sequence valued map $\Omega\ni z\mapsto \left(w^k_0(z),w^k_{-1}(z),w^k_{-2}(z),\dots\right)$ ($0 \leq k \leq 2$) and $\mathscr{F}=\left(\mathscr{F}_0,\mathscr{F}_1,\mathscr{F}_2,0,0,0,\dots\right)$,
\begin{subequations}\label{eq;system in terms of fourier coefficients}
    \begin{align}
    &\overline{\partial}w_{1}^{0}(z) + \partial w_{-1}^{0}(z) = \mathscr{F}_{0}(z),\label{eq; w_{-1}^{0} relation btw components of the function F}\\
    &\overline{\partial} w_{0}^{0}(z) + \partial w_{-2}^{0}(z) = \mathscr{F}_{1}(z),\\
    &\overline{\partial} w_{-1}^{0}(z) + \partial w_{-3}^{0}(z) = \mathscr{F}_{2}(z),\\
    &\overline{\partial} w_{-n}^{0}(z) + \partial w_{-(n+2)}^{0}(z) = 0, \quad n\geq 2,\\
    &\overline{\partial} w_{-n}^{1}(z) + \partial w_{-(n+2)}^{1}(z) = w_{-(n+1)}^{0}(z), \quad  n\in \mathbb{Z},\label{eq; w_{-1}^{0} relation btw components of the function}\\
    &\overline{\partial} w_{-n}^{2}(z) + \partial w_{-(n+2)}^{2}(z) = w_{-(n+1)}^{1}(z), \quad  n\in \mathbb{Z},\label{eq; v_0^{-1} relation btw components of the function}
    \end{align}  
\end{subequations}
subject to
    \begin{align}
         w_{-n}^{k}|_{\Gamma} = g_{-n}^{k} \quad (n\in \mathbb{Z},\; 0\leq k \leq 2),
    \end{align}
    where $\vg^k = \vw^k|_{\Gamma}$ $(0\leq k \leq 2)$ and $\vg = \left(g_{0}, g_{-1}, g_{-2},...\right)$ is the sequence with non-positive indices of Fourier coefficients of $\vg$ given by \begin{align}
        g_{-n}(z) = \frac{1}{2\pi} \int\limits_{0}^{2\pi} \vg(z,\vuu_{\varphi}) \mathrm{e}^{\mathrm{i}n\varphi} \, d\varphi \quad  (n\geq 0). 
    \end{align}
Rewriting the system of equation \eqref{eq; w_{-1}^{0} relation btw components of the function F}-\eqref{eq; v_0^{-1} relation btw components of the function} in terms of the left shift operator, we have
    \begin{subequations}\label{BVP for non-attenuated case}
        \begin{align}
            \overline{\partial} \overline{w_{-1}^{0}}(z) + \partial w_{-1}^{0}(z) &= \mathscr{F}_{0}(z),\\
            \overline{\partial} \vw^{0}(z) + \mathcal{L}^2\partial \vw^{0}(z) &= \mathcal{L}\mathscr{F}(z),\\
            \overline{\partial} \vw^{1}(z) + \mathcal{L}^2\partial \vw^{1}(z) &= \mathcal{L}\vw^{0}(z),\label{BVP for non-attenuated case, 32(c)}\\
            \overline{\partial} \vw^{2}(z) + \mathcal{L}^2\partial \vw^{2}(z) &= \mathcal{L}\vw^{1}(z),\label{BVP for non-attenuated case, 32(d)}
        \end{align}
    \end{subequations}
    subject to
    \begin{align}
          \vw^k|_{\Gamma}= \vg^k \quad (0\leq k \leq 2).
    \end{align}
Note that $\mathcal{L}^3\mathscr{F} = \textbf{0} = (0,0,0,...)$. So, $\mathcal{L}^2 \vw^{0}$ can be obtained by solving the following BVP using equations \eqref{eq;homo Beltrami like eq} and \eqref{Bukhgeim-Beltrami operator},
     \begin{align}\label{BVP in w^0}
         \begin{rcases}
              \overline{\partial} \mathcal{L}^2\vw^{0}(z) + \mathcal{L}^2\partial \mathcal{L}^2 \vw^{0}(z) = \textbf{0}\\
              \mathcal{L}^2 \vw^{0}|_{\Gamma} =  \mathcal{L}^2 \vg^{0},
                 \end{rcases}
     \end{align}
and then knowing $\mathcal{L}^2 \vw^{0}$, $\mathcal{L}^1 \vw^{1}$ can be obtained by solving the following BVP using equation \eqref{eq;nonhom_} and Proposition \ref{prop_non-homogeneous solution},
     \begin{align}\label{BVP in w^1}
         \begin{rcases}
              \overline{\partial} \mathcal{L}^1\vw^{1}(z) + \mathcal{L}^2\partial \mathcal{L}^1 \vw^{1}(z) = \mathcal{L}^2 \vw^{0}(z)\\
              \mathcal{L}^1 \vw^{1}|_{\Gamma} =  \mathcal{L}^1 \vg^{1}.
                \end{rcases}
     \end{align}
Again as $\mathcal{L}^1 \vw^{1}$ is known, using equation \eqref{eq;nonhom_} and Proposition \ref{prop_non-homogeneous solution}, $\vw^{2}$ is obtained by solving the following BVP
     \begin{align}\label{BVP in w^2}
         \begin{rcases}
              \overline{\partial} \vw^{2}(z) + \mathcal{L}^2\partial  \vw^{2}(z) = \mathcal{L}^1 \vw^{1}(z)\\
               \vw^{2}|_{\Gamma} =   \vg^{2}
               \end{rcases}
     \end{align}
Therefore, from equations \eqref{BVP in w^0}, \eqref{BVP in w^1} and \eqref{BVP in w^2}, the solution of the above BVP \eqref{BVP for non-attenuated case} satisfies
   \begin{align}
       \mathcal{L}^{2-k}\vw^k(z) = \sum\limits_{j=0}^{k} \mathcal{T}^j\mathcal{L}^{2-k+j}\left[ \mathfrak{B}\vg^{k-j}\right](z) \quad (z\in \Omega,~ 0\leq k \leq 2),
   \end{align}
    with the following estimate which follows from \cite[Prop. 4.1]{fujiwara2024inversion}
    \begin{align}\label{estimate btwn w^k and g^j}
        \left\lVert \mathcal{L}^{2-k}\vw^k\right\rVert_{2-k,k+1}^2 \lesssim \sum\limits_{j=0}^{k} \left\lVert \mathcal{L}^{2-j}\vg^j\right\rVert_{\frac{7}{2},j+\frac{1}{2}}^2 \quad (0\leq k \leq 2).
    \end{align}
Now, it remains to find $w_{0}^{0}$,  $w_{-1}^{0}$ and $w_{0}^{1}$. Taking $n=-1,0$ in the equation \eqref{eq; w_{-1}^{0} relation btw components of the function} and $n=-1$ in the equation \eqref{eq; v_0^{-1} relation btw components of the function}, we get \begin{align}\label{eq; value of w_0^0}
    w_{0}^{0}(z) = \overline{\partial} \overline{w_{-1}^{1}}(z) + \partial  w_{-1}^{1}(z), \;\; w_{-1}^{0}(z) = \overline{\partial} w_0^{1}(z) + \partial  w_{-2}^{1}(z) \;\; \text{and} \;\; w_0^{1}(z) = \overline{\partial} \overline{w_{-1}^{2}}(z) + \partial  w_{-1}^{2}(z).
\end{align}
Thus, all the values of the sequences $\vw^{0}, \vw^{1}$ and $\vw^{2}$ are now known and hence from the equation \eqref{eq; w_{-1}^{0} relation btw components of the function F}--\eqref{eq; v_0^{-1} relation btw components of the function} the values of $\mathscr{F}_{0}$, $\mathscr{F}_{1}$ and $\mathscr{F}_{2}$ can be recovered. Finally, we can uniquely recover the vector field ${\vf}$ and symmetric 2-tensor field ${\vF}$ by
\begin{align}
    {\vf} := (2 \mathbb{R}e\mathscr{F}_1, 2 \mathbb{I}m\mathscr{F}_1) \quad \text{and} \quad {\vF} := \begin{bmatrix}
        \mathscr{F}_{0} + 2 \mathbb{R}e\mathscr{F}_{2} & 2 \mathbb{I}m\mathscr{F}_2\\
        2 \mathbb{I}m\mathscr{F}_2 & \mathscr{F}_{0} - 2 \mathbb{R}e\mathscr{F}_{2}
    \end{bmatrix},
\end{align}
where
\begin{align}\label{estimate btwn w^0 and w^2}
     \mathscr{F}_{0} = 2 \mathbb{R}e\left[ \partial w_{-1}^{0}\right], \quad \mathscr{F}_{1} = \overline{\partial} w_{0}^{0} + \partial w_{-2}^{0} \quad \text{and} \quad \mathscr{F}_{2} = \overline{\partial} w_{-1}^{0} + \partial w_{-3}^{0}.
\end{align}
So, we have recovered $\vf$ and $\vF$ from the set of AMRTs. This concludes the recovery part. Next we prove the stability estimate. Using the expression $\mathscr{F} = (\mathscr{F}_0,\mathscr{F}_1,\mathscr{F}_3,0,0,\cdots) $, equation \eqref{equation;relation btw components of the function} and the fact that each component of $\vf$ and $\vF$ are real valued, we have  
\begin{align*}
\left\lVert\mathscr{F}\right\rVert_{L^{2}(\Omega)}^2 &\vcentcolon= \lVert (\mathscr{F}_0,\mathscr{F}_1,\mathscr{F}_3,0,0,\dots)\rVert_{L^{2}(\Omega)}^2 \\ 
&\vcentcolon= \left\lVert\mathscr{F}_{0}\right\rVert_{L^{2}(\Omega)}^2 + \left\lVert\mathscr{F}_{1}\right\rVert_{L^{2}(\Omega)}^2 + \left\lVert\mathscr{F}_{2}\right\rVert_{L^{2}(\Omega)}^2\\
&= \left\lVert\frac{{F}_{11} + {F}_{22}}{2}\right\rVert_{L^{2}(\Omega)}^2  + \left\lVert\frac{{f}_{1}+ \mathrm{i}{f}_{2}}{2}\right\rVert_{L^{2}(\Omega)}^2 + \left\lVert\frac{{F}_{11} - {F}_{22}}{4} + \mathrm{i} \frac{{F}_{12}}{2}\right\rVert_{L^{2}(\Omega)}^2\\
& \geq \frac{1}{16}\left( \left\lVert{F}_{11} + {F}_{22}\right\rVert_{L^{2}(\Omega)}^2 + \left\lVert{F}_{11} - {F}_{22}\right\rVert_{L^{2}(\Omega)}^2\right) +\frac{1}{4}\left( \left\lVert{F}_{12}\right\rVert_{L^{2}(\Omega)}^2 + \left\lVert{f}_{1}\right\rVert_{L^{2}(\Omega)}^2 + \left\lVert{f}_{2}\right\rVert_{L^{2}(\Omega)}^2\right).
\end{align*}
Using the parallelogram law, we get
\begin{align}\label{relation mathscr{F} and vetor field and 2-tensor field}
\left\lVert\mathscr{F}\right\rVert_{L^{2}(\Omega)}^2
 & \geq \frac{1}{8}\left( \left\lVert{F}_{11} \right\rVert_{L^{2}(\Omega)}^2 + \left\lVert {F}_{22}\right\rVert_{L^{2}(\Omega)}^2\right) +\frac{1}{4}\left( \left\lVert{F}_{12}\right\rVert_{L^{2}(\Omega)}^2 + \left\lVert{f}_{1}\right\rVert_{L^{2}(\Omega)}^2 + \left\lVert{f}_{2}\right\rVert_{L^{2}(\Omega)}^2\right)\nonumber\\
 &\gtrsim \left( \left\lVert{F}_{11} \right\rVert_{L^{2}(\Omega)}^2 + \left\lVert {F}_{22}\right\rVert_{L^{2}(\Omega)}^2 +2 \left\lVert{F}_{12}\right\rVert_{L^{2}(\Omega)}^2\right) + \left(\left\lVert{f}_{1}\right\rVert_{L^{2}(\Omega)}^2 + \left\lVert{f}_{2}\right\rVert_{L^{2}(\Omega)}^2\right) \nonumber\\
&= \left\lVert{\vF}\right\rVert_{L^{2}(\Omega)}^2 + \left\lVert{\vf }\right\rVert_{L^{2}(\Omega)}^2
\end{align}
Thus, we have $ \left\lVert{\vF}\right\rVert_{L^{2}(\Omega)}^2 + \left\lVert{\vf }\right\rVert_{L^{2}(\Omega)}^2 \lesssim \left\lVert\mathscr{F}\right\rVert_{L^{2}(\Omega)}^2$. So, it suffices to estimate $\mathscr{F}$ by the given data.
Observe that repeated differentiation of \eqref{BVP for non-attenuated case, 32(c)} and \eqref{BVP for non-attenuated case, 32(d)} gives $$\nabla^q (\mathcal{L}\vw^{k-1})= \overline{\partial} [\nabla^q \vw^k] + \mathcal{L}^2 \partial[\nabla^q \vw^k] \quad (q\in\mathbb{N}),$$ where $\nabla$ is the vector differential operator and it give the estimates $\left\lVert \mathcal{L}\vw^{k-1}\right\rVert_{0,q}^2 \lesssim  \left\lVert \vw^{k}\right\rVert_{0,q+1}^2$ for $k=1,2$ and $q\in\mathbb{N}$, and using it along with \eqref{eq; value of w_0^0}, we have 
\begin{align} \label{eq: xxx}
     &\left\lVert \vw^0\right\rVert_{0,1}^2 = \lVert w_0^0\rVert^2_{H^1} + \sum_{n\in\mathbb{N}} \lVert w_{-n}^0\rVert^2_{H^1} = \lVert \overline{\partial} \overline{w^1_{-1}} + \partial w^1_{-1}\rVert^2_{H^1} + \lVert \mathcal{L}\vw^0 \rVert^2_{0,1} \lesssim \lVert \vw^1\rVert^2_{0,2} \nonumber \\  
     \text{and} \quad &\lVert \vw^1\rVert^2_{0,2} = \lVert w_0^1\rVert^2_{H^2} + \sum_{n\in\mathbb{N}} \lVert w_{-n}^1\rVert^2_{H^2} = \lVert \overline{\partial} \overline{w^2_{-1}} + \partial w^2_{-1}\rVert^2_{H^2} + \lVert \mathcal{L}\vw^1 \rVert^2_{0,1} \lesssim  \left\lVert \vw^2\right\rVert_{0,3}^2 \nonumber \\ 
     \text{giving} \quad &\left\lVert \vw^0\right\rVert_{0,1} \lesssim  \left\lVert \vw^2\right\rVert_{0,3}.
\end{align}
From equations \eqref{estimate btwn w^k and g^j}, \eqref{estimate btwn w^0 and w^2} and \eqref{eq: xxx}, we have
\begin{align}\label{relation mathscr{F} and components of vetor field and 2-tensor field}
\left\lVert\mathscr{F}\right\rVert_{0,0}^2 \lesssim \left\lVert \vw^{0}\right\rVert_{0,1}^2 \lesssim \left\lVert \vw^{2}\right\rVert_{0,3}^2 \lesssim \sum\limits_{j=0}^{2} \left\lVert \mathcal{L}^{2-j}\vg^{j}\right\rVert_{\frac{7}{2},j+\frac{1}{2}},
\end{align}
and by \eqref{system of transport eq.}, we have 
\begin{align}
\sum\limits_{j=0}^{2} \left\lVert \mathcal{L}^{2-j}\vg^{j}\right\rVert_{\frac{7}{2},j+\frac{1}{2}} \lesssim\sum\limits_{k=0}^{2} \left\lVert \mathcal{M}^{(k)}\mathscr{F}\right\rVert_{\frac{7}{2},k+\frac{1}{2}}.    
\end{align}
Finally, by the equations \eqref{relation mathscr{F} and vetor field and 2-tensor field} and \eqref{relation mathscr{F} and components of vetor field and 2-tensor field},
\begin{align}
    \left\lVert{\vf }\right\rVert_{L^{2}(\Omega)}^2  + \left\lVert{\vF}\right\rVert_{L^{2}(\Omega)}^2  \lesssim \sum\limits_{k=0}^{2} \left\lVert \mathcal{M}^{(k)}\mathscr{F}\right\rVert_{\frac{7}{2},k+\frac{1}{2}}.
\end{align}
This completes the proof for the non-attenuated case.

\subsection{The attenuated case: }\label{sec; attenuated case}
As mentioned earlier, for the attenuated case, we assume that $\mathfrak{a}>0$ in $\overline{\Omega}$. First, as in \cite{sadiq2015range}, we introduce a special integrating factor function $\mathscr{H}$, which helps to define a one-to-one correspondence between the $\mathcal{L}^2$-analytic map $\vw := (w_{0}, w_{-1}, w_{-2},\dots)$ that satisfies \eqref{eq;homo Beltrami like eq} and the solution $\vv := (v_{0}, v_{-1}, v_{-2},\dots)$ of the equation $\overline{\partial}\vv + \mathcal{L}^2\partial \vv + \mathfrak{a}\mathcal{L}\vv = 0$. The function $\mathscr{H}$, appeared first time in \cite{natterer2001mathematics}, is defined as 
\begin{align}
    \mathscr{H}(z, \vuu_{\varphi}) \vcentcolon= \int\limits_{0}^{\infty}\mathfrak{a}(z+t\vuu_{\varphi})\,dt - \frac{1}{2}(I - \mathrm{i}\mathcal{H})\mathcal{R}\mathfrak{a}\left(z\cdot{\vuu_{\varphi}}^{\perp}, {\vuu_{\varphi}}^{\perp}\right),
\end{align}
where the notations $\mathcal{H}$ and $\mathcal{R}$ stands for the Hilbert transform and the Radon transform, respectively, and given by
\begin{align}
    \mathcal{H}f(s, \vuu_{\varphi}) = \frac{1}{\pi} \int\limits_{-\infty}^{\infty} \frac{f (t,\vuu_{\varphi})}{s-t}\,dt \quad \text{and} \quad \mathcal{R}\mathfrak{a}(s, \vuu_{\varphi}) = \int\limits_{-\infty}^{\infty} \mathfrak{a}\left(s{\vuu_{\varphi}}^{\perp}+ t{\vuu_{\varphi}}\right)\,dt.
\end{align}

From \cite{finch2003attenuated,natterer2001mathematics}, it follows that the negative Fourier coefficients for the function $\mathscr{H}$ are vanishing and have the expansion 
\begin{align*}
    \mathrm{e}^{-\mathscr{H}(z,\vuu_{\varphi})} \vcentcolon= \sum\limits_{n=0}^{\infty}\alpha_{n}(z)\mathrm{e}^{\mathrm{i}n\varphi},\quad \mathrm{e}^{\mathscr{H}(z,\vuu_{\varphi})} \vcentcolon= \sum\limits_{n=0}^{\infty}\beta_{n}(z)\mathrm{e}^{\mathrm{i}n\varphi} \quad ((z,\vuu_{\varphi})\in \overline{\Omega}\times\mathbb{S}^1).
\end{align*}
Define the sequence-valued maps using the Fourier coefficients of the functions $e^{\pm\mathscr{H}}$ by
\begin{align*}
 z \mapsto \alpha(z) \vcentcolon= \left(\alpha_{0}(z), \alpha_{1}(z), \alpha_{2}(z),\dots\right), \quad   z \mapsto \beta(z) \vcentcolon= \left(\beta_{0}(z), \beta_{1}(z), \beta_{2}(z),\dots\right), 
\end{align*}
and the operator $\mathrm{e}^{\pm \mathcal{G}}$, as in \cite{fujiwara2019fourier}, component-wise for each $n\leq0$ by 
\begin{align}\label{eq;operators e}
    (\mathrm{e}^{-\mathcal{G}}\vv)_{n} = (\alpha \ast \vv)_n = \sum\limits_{m=0}^{\infty}\alpha_mv_{n-m} \quad\text{and}\quad 
    (\mathrm{e}^{\mathcal{G}}\vv)_{n} = (\beta \ast \vv)_n = \sum\limits_{m=0}^{\infty}\beta_mv_{n-m}.
    \end{align}
We also have the commutating property $[\mathrm{e}^{\mathcal{\pm G}}, \mathcal{L}] = 0.$

Next, we state a lemma which establishes a relation between attenuated and non-attenuated cases. 

\begin{lemma}\label{conerter result} \cite[Lemma 5.1]{fujiwara2024inversion}
Let $\mathfrak{a} \in C^{1, \alpha}(\overline{\Omega})$, for $1/2 <\alpha \leq 1$, and let $\mathrm{e}^{\pm \mathcal{G}}$ be the operators as defined in \eqref{eq;operators e}.
\begin{enumerate}
    \item[(i)] If $\vv\in C^{1}(\Omega, \ell^1)$ solves $\overline{\partial}\vv + \mathcal{L}^2\partial \vv + \mathfrak{a}\mathcal{L}\vv = \vv$, then $\vw = \mathrm{e}^{-\mathcal{G}}\vv \in C^{1}(\Omega, \ell^1) $ solves $\overline{\partial}\vw + \mathcal{L}^2\partial \vw = \mathrm{e}^{-\mathcal{G}}\vv$.
    
    \item[(ii)] Conversely, if $\vw \in C^{1}(\Omega,\ell^1) $ solves $\overline{\partial}\vw + \mathcal{L}^2\partial \vw = \mathrm{e}^{-\mathcal{G}}\vv$, then $\vv = \mathrm{e}^{\mathcal{G}} \vw\in C^{1}(\Omega,\ell^1)$ solves $\overline{\partial}\vv + \mathcal{L}^2\partial \vv + \mathfrak{a}\mathcal{L}\vv = \vv$.
\end{enumerate}
   
\end{lemma}
 
As done in the non-attenuated case, we rewrite the system of equations in \eqref{eq; w_{-1}^{0} relation btw components of the function F}-\eqref{eq; v_0^{-1} relation btw components of the function} in terms of the left shift operator as
    \begin{subequations}
        \begin{align}
            \overline{\partial} \overline{v_{-1}^{0}}(z) + \partial v_{-1}^{0}(z) + \mathfrak{a} v_{0}^{0}(z) &= \mathscr{F}_{0}(z),\\
              \overline{\partial} \vv^{0}(z) + \mathcal{L}^2\partial \vv^{0}(z) + \mathfrak{a}\mathcal{L}\vv^{0}(z) &= \mathcal{L}\mathscr{F}(z),\\
              \overline{\partial} \vv^{1}(z) + \mathcal{L}^2\partial \vv^{1}(z) + \mathfrak{a}\mathcal{L}\vv^{1}(z) &= \mathcal{L}\vv^{0}(z),\\
              \overline{\partial} \vv^{2}(z) + \mathcal{L}^2\partial \vv^{2}(z) + \mathfrak{a}\mathcal{L}\vv^{2}(z) &= \mathcal{L}\vv^{1}(z),
        \end{align}
    \end{subequations}
        subject to
    \begin{align}
         \vg^k = \vv^k|_{\Gamma} \quad (0\leq k \leq 2).
    \end{align}
By equation \eqref{eq;operators e} and Lemma \ref{conerter result}, it follows that for $\vw = \mathrm{e}^{-\mathcal{G}}\vv$, we have
\begin{subequations}
    \begin{align}
        \overline{\partial} \overline{w_{-1}^{0}}(z) + \partial w_{-1}^{0}(z) &= (\mathrm{e}^{-\mathcal{G}}\mathscr{F})_{0}(z),\\
        \overline{\partial} \vw^{0}(z) + \mathcal{L}^2\partial \vw^{0}(z) &= \mathcal{L}[\mathrm{e}^{-\mathcal{G}}\mathscr{F}](z),\\
        \overline{\partial} \vw^{1}(z) + \mathcal{L}^2\partial \vw^{1}(z) &= \mathcal{L}\vw^{0}(z),\\
        \overline{\partial} \vw^{2}(z) + \mathcal{L}^2\partial \vw^{2}(z) &= \mathcal{L}\vw^{1}(z),
    \end{align}
\end{subequations}
subject to
\begin{align}
       \vw^k|_{\Gamma} =  \mathrm{e}^{-\mathcal{G}}\vg^k, \quad 0\leq k \leq 2.
\end{align}
Again from the fact $\mathcal{L}^3\mathscr{F} = \textbf{0} = (0,0,0,...)$ and that $\mathrm{e}^{\pm\mathcal{G}}$ commute with the operator $\mathcal{L}$, we have
    \begin{align}
        \mathcal{L}^3[\mathrm{e}^{-\mathcal{G}}\mathscr{F}] = \mathrm{e}^{-\mathcal{G}}\mathcal{L}^3\mathscr{F} = \mathrm{e}^{-\mathcal{G}}\textbf{0} = \textbf{0}.
    \end{align}
Since the attenuation coefficient $\mathfrak{a}$ and $\vg^k{|_\Gamma}$ are known for $0\leq k \leq 2$, from the equations \eqref{eq;traces} and \eqref{eq; solution of homo Beltrami like eq}, it follows that $\mathcal{L}^2 \vw^{0}$ can be obtained by solving the following BVP 
\begin{align}
    \begin{rcases}
        \overline{\partial} \mathcal{L}^2\vw^{0}(z) + \mathcal{L}^2\partial \mathcal{L}^2 \vw^{0}(z) &= \textbf{0}\\
        \mathcal{L}^2 \vw^{0}|_{\Gamma} =  \mathrm{e}^{-\mathcal{G}}\mathcal{L}^2 \vg^{0}
    \end{rcases},
\end{align}
$\mathcal{L}^1 \vw^{1}$ can be obtained by solving the following BVP 
     \begin{align}
         \begin{rcases}
              \overline{\partial} \mathcal{L}^1\vw^{1}(z) + \mathcal{L}^2\partial \mathcal{L}^1 \vw^{1}(z) &= \mathcal{L}^2 \vw^{0}(z)\\
              \mathcal{L}^1 \vw^{1}|_{\Gamma} =  \mathrm{e}^{-\mathcal{G}}\mathcal{L}^1 \vg^{1},
        \end{rcases}
     \end{align}
and $ \vw^{2}$ can be obtained by solving the following BVP 
     \begin{align}
        \begin{rcases}
            \overline{\partial} \vw^{2}(z) + \mathcal{L}^2\partial  \vw^{2}(z) &= \mathcal{L}^1 \vw^{1}(z)\\
            \vw^{2}|_{\Gamma} =  \mathrm{e}^{-\mathcal{G}} \vg^{2}.
        \end{rcases}
     \end{align}
So, for each $0\leq k \leq 2$, $\mathcal{L}^{2-k}\vw$ are now known, and using techniques as in the non-attenuated case, we get the values of $\vw^{k}$. Then, by Lemma \ref{conerter result}, we get $\vv^k = \mathrm{e}^{\mathcal{G}}\vw^k$ and $\mathscr{F}_{0}, \mathscr{F}_{1}, \mathscr{F}_{2}$ where $\mathscr{F} = \mathrm{e}^{\mathcal{G}}[\mathrm{e}^{-\mathcal{G}}\mathscr{F}] $.

Hence the vector field ${\vf}$ and the symmetric 2-tensor field ${\vF}$ can be obtained by the following
\begin{align}
    {\vf} := \left(2 \mathbb{R}e\mathscr{F}_1, 2 \mathbb{I}m\mathscr{F}_1\right) \quad \text{and} \quad
    {\vF}(x) := \begin{bmatrix}
        \mathscr{F}_{0} + 2 \mathbb{R}e\mathscr{F}_{2} & 2 \mathbb{I}m\mathscr{F}_2\\
        2 \mathbb{I}m\mathscr{F}_2 & \mathscr{F}_{0} - 2 \mathbb{R}e\mathscr{F}_{2}
    \end{bmatrix},
\end{align}
where
\begin{align}
    (\mathrm{e}^{-\mathcal{G}}\mathscr{F})_{0}:= 2 \mathbb{R}e\left[ \partial w_{-1}^{0}\right], \;\; (\mathrm{e}^{-\mathcal{G}}\mathscr{F})_{2}(x) = \overline{\partial} w_{-1}^{0} + \partial w_{-3}^{0} \;\; \text{and} \;\; 
    (\mathrm{e}^{-\mathcal{G}}\mathscr{F})_{1}(x) = \overline{\partial} w_{0}^{0} + \partial w_{-2}^{0},
\end{align}
with the estimate 
\begin{align}\label{relation mathscr{F} and components of vetor field and 2-tensor field for the attenuated case}
\left\lVert e^{-\mathcal{G}}\mathscr{F}\right\rVert_{0,0}^2  \lesssim \sum\limits_{j=0}^{2} \left\lVert e^{-\mathcal{G}}\mathcal{L}^{2-j}g^{j}\right\rVert_{\frac{7}{2},j+\frac{1}{2}}.
\end{align}
Using the fact that $\mathscr{F} = e^{\mathcal{G}}\left[e^{-\mathcal{G}} \mathscr{F}\right]$, and Lemma \ref{conerter result}, we get the same estimate as we got for the non-attenuated case, $\left\lVert \mathscr{F}\right\rVert_{0,0}^2  \lesssim \sum\limits_{j=0}^{2} \left\lVert \mathcal{L}^{2-j}g^{j}\right\rVert_{\frac{7}{2},j+\frac{1}{2}}$. This completes the proof of the main theorem.

\subsection{Proof of the proposition \ref{proposition; transport model}} \label{sec,proof of Prop}
To conclude and make the article self-contained, we now prove Proposition \ref{proposition; transport model}.

For $(\vx,\vuu_{\varphi})\in \Omega\times \mathbb{S}^{1}$ and $k=1,2$, equation \eqref{eq:system of BVP} gives
\begin{subequations}
    \begin{align}
    &\frac{d}{dt}\left[ \mathrm{e}^{-\int\limits_{t}^{\infty}\mathfrak{a}(\vx +s\vuu_{\varphi})\,ds}\vv^{0}(\vx +t\vuu_{\varphi},\vuu_{\varphi})\right] = \mathrm{e}^{-\int\limits_{t}^{\infty}\mathfrak{a}(\vx +s\vuu_{\varphi})\,ds} \left({\vf}(\vx +t\vuu_\varphi)\cdot\vuu_{\varphi}+\left\langle {\vF} (\vx +t\vuu_\varphi), \vuu_\varphi^{2} \right\rangle\right)\label{eq: derivative of v^0}\\ 
    &\text{and} \quad  \frac{d}{dt}\left[ \mathrm{e}^{-\int\limits_{t}^{\infty}\mathfrak{a}(\vx +s\vuu_{\varphi})\,ds}\vv^{k}(\vx +t\vuu_{\varphi},\vuu_{\varphi})\right] = \mathrm{e}^{-\int\limits_{t}^{\infty}\mathfrak{a}(\vx +s\vuu_{\varphi})\,ds}\vv^{k-1}(\vx +t\vuu_{\varphi},\vuu_{\varphi})\label{eq: derivative of v^k}.
\end{align}
\end{subequations}
Now, we integrate both the equations \eqref{eq: derivative of v^0} and \eqref{eq: derivative of v^k}, along the line which passes through $\vx$ in the direction of $\vuu_\varphi$, we get
\begin{align} \label{eqpp0}
    \mathrm{e}^{-\int\limits_{\vx\cdot\vuu_{\varphi}}^{\infty}\mathfrak{a}(\Pi_{\varphi}(\vx) +s\vuu_{\varphi})\,ds}\vv^{0}(\vx +t\vuu_{\varphi},\vuu_{\varphi}) 
    &= \int\limits_{-\infty}^{x\cdot\vuu_{\varphi}}\frac{d}{dt}\left[ \mathrm{e}^{-\int\limits_{t}^{\infty}\mathfrak{a}(\Pi_{\varphi}(\vx) +s\vuu_{\varphi})\,ds}\vv^{0}(\Pi_{\varphi}(\vx)+t\vuu_{\varphi},\vuu_{\varphi})\right]\,dt \nonumber \\
    &= \int\limits_{-\infty}^{\vx\cdot\vuu_{\varphi}}\mathrm{e}^{-\int\limits_{t}^{\infty}\mathfrak{a}(\Pi_{\varphi}(\vx) +s\vuu_{\varphi})\,ds} \mathcal{F}(\Pi_{\varphi}(\vx) +t\vuu_{\varphi})\,dt,
\end{align}
 where $\mathcal{F}(\Pi_{\varphi}(\vx) +t\vuu_{\varphi})  = {\vf}(\Pi_{\varphi}(\vx) +t\vuu_\varphi)\cdot\vuu_{\varphi}+\left\langle {\vF} (\Pi_{\varphi}(\vx) +t\vuu_\varphi), \vuu_\varphi^{2} \right\rangle.$\\
 Here, observe that $\int\limits_{\vx\cdot\vuu_\varphi}^{\infty}\mathfrak{a}(\Pi_{\varphi}(\vx) +s\vuu_{\varphi})\,ds = \int\limits_{0}^{\infty}\mathfrak{a}(\vx + s\vuu_{\varphi})\,ds  $ 
and using the integration by parts formula in equation \eqref{eq:system of BVP} gives
\begin{align}
     & \mathrm{e}^{-\int\limits_{\vx\cdot\vuu_{\varphi}}^{\infty}\mathfrak{a}(\Pi_{\varphi}(\vx) +s\vuu_{\varphi})\,ds}\vv^{1}(\vx +t\vuu_{\varphi},\vuu_{\varphi}) \nonumber 
    = \int\limits_{-\infty}^{\vx\cdot\vuu_{\varphi}}\frac{d}{dt}\left[ \mathrm{e}^{-\int\limits_{t}^{\infty}\mathfrak{a}(\Pi_{\varphi}(\vx) +s\vuu_{\varphi})\,ds}\vv^{1}(\Pi_{\varphi}(\vx)+t\vuu_{\varphi},\vuu_{\varphi})\right]\,dt \nonumber \\
    &= \int\limits_{-\infty}^{\vx\cdot\vuu_{\varphi}}\mathrm{e}^{-\int\limits_{t}^{\infty}\mathfrak{a}(\Pi_{\varphi}(\vx) +s\vuu_{\varphi})\,ds}\vv^{0}(\Pi_{\varphi}(\vx) +t\vuu_{\varphi})\,dt \nonumber\\
    &= \mathrm{e}^{-\int\limits_{\vx\cdot\vuu_{\varphi}}^{\infty}\mathfrak{a}(\Pi_{\varphi}(\vx) +s\vuu_{\varphi})\,ds}(\vx\cdot\vuu_{\varphi})\vv^{0}(\Pi_{\varphi}(\vx) +t\vuu_{\varphi}) - \int\limits_{-\infty}^{\vx\cdot\vuu_{\varphi}}t\mathrm{e}^{-\int\limits_{t}^{\infty}\mathfrak{a}(\Pi_{\varphi}(\vx) +s\vuu_{\varphi})\,ds}\mathcal{F}(\Pi_{\varphi}(\vx) +t\vuu_{\varphi})\,dt
\end{align}
and 
\begin{align}
    &\mathrm{e}^{-\int\limits_{\vx\cdot\vuu_{\varphi}}^{\infty}\mathfrak{a}(\Pi_{\varphi}(\vx) +s\vuu_{\varphi})\,ds}  \vv^{2}(\vx +t\vuu_{\varphi},\vuu_{\varphi}) \\
    =\ &  \mathrm{e}^{-\int\limits_{\vx\cdot\vuu_{\varphi}}^{\infty}\mathfrak{a}(\Pi_{\varphi}(\vx) +s\vuu_{\varphi})\,ds}\sum\limits_{n=1}^{2} \frac{(-1)^{n-1}}{n!} (\vx\cdot\vuu_{\varphi})^{n} \vv^{2-n}(\Pi_{\varphi}(\vx) +t\vuu_{\varphi})\nonumber\\ 
    &+ \int\limits_{-\infty}^{\vx\cdot\vuu_{\varphi}}\frac{t^{2}}{2}\mathrm{e}^{-\int\limits_{t}^{\infty}\mathfrak{a}(\Pi_{\varphi}(\vx) +s\vuu_{\varphi})\,ds}\mathcal{F}(\Pi_{\varphi}(\vx) +t\vuu_{\varphi})\,dt.
\end{align}
For $k\in\{1,2\}$, the above equations gives
\begin{align} \label{eqpp1}
    \vv^{k}(\vx +t\vuu_{\varphi},\vuu_{\varphi}) \nonumber
    &= \sum\limits_{n=1}^{k} (-1)^{n-1}\frac{(\vx\cdot\vuu_{\varphi})^{n}}{n!}\vv^{k-n}(\Pi_{\varphi}(\vx) +t\vuu_{\varphi})\nonumber\\ & ~~~~ +(-1)^{k} \int\limits_{-\infty}^{\vx\cdot\vuu_{\varphi}}\frac{t^{k}}{k!} \mathrm{e}^{-\int\limits_{t}^{\vx \cdot \vuu_\varphi} \mathfrak{a}(\Pi_{\varphi}(\vx) +s\vuu_{\varphi})\,ds}\mathcal{F}(\Pi_{\varphi}(\vx) +t\vuu_{\varphi})\,dt.
\end{align}
If $(\vx,\vuu_\varphi) \in \Gamma_+$ and $t>\vx\cdot\vuu_\varphi$, then $\vf(\vx+(t-\vx\cdot\vuu_\varphi)\vuu_\varphi)=0=\vF(\vx+(t-x\cdot\vuu_\varphi)\vuu_\varphi)$ and so
\begin{align} \label{eqpp2}
   \mathcal{M}_a^{(k)}\mathcal{F}(\vx,\vuu_\varphi) &= \int\limits_{-\infty}^{\infty}{t^{k}} \mathrm{e}^{-\int\limits_{t}^{\vx \cdot \vuu_\varphi} \mathfrak{a}(\Pi_{\varphi}(\vx) +s\vuu_{\varphi})\,ds}\mathcal{F}(\Pi_{\varphi}(\vx) +t\vuu_{\varphi})\,dt \nonumber
   \\ &=\int\limits_{-\infty}^{\vx\cdot\vuu_{\varphi}}{t^{k}} \mathrm{e}^{-\int\limits_{t}^{\vx \cdot \vuu_\varphi} \mathfrak{a}(\Pi_{\varphi}(\vx) +s\vuu_{\varphi})\,ds}\mathcal{F}(\Pi_{\varphi}(\vx) +t\vuu_{\varphi})\,dt.
\end{align}
The set of equations \eqref{eqpp0}--\eqref{eqpp2} gives \eqref{system of transport eq.}. Also, if $s\geq1$, then ${\vf} \in H_0^{s}(\Omega; S^1)$, ${\vF} \in H_0^{s}(\Omega; S^2)$ and $\mathfrak{a} \in C^{1,\alpha}(\overline{\Omega})$, $1/2<\alpha\leq1$ along with equations \eqref{eq:system of BVP}, \eqref{eqpp0} and \eqref{eqpp1} implies that $\vv^k \in H^s(\Omega\times\mathbb{S}^1)$ for $0\leq k \leq 2$.


\section{Acknowledgements}\label{sec:acknowledgements}
The first author gratefully acknowledges the Senior Research Fellowship from the UGC, the Government of India. The second author gratefully acknowledges the Post Doctoral Fellowship under the ISIRD project 9--551/2023/IITRPR/10229 from IIT Ropar. The authors are grateful to Manmohan Vashisth for his careful reading and valuable suggestions which enhanced the quality of the article. Finally, first author would like to thank Venky Krishnan and Rohit Kumar Mishra for several interesting discussions regarding stability estimate for the aforementioned ray transform.

\addcontentsline{toc}{section}{References}
\bibliography{References}

\newcommand{\etalchar}[1]{$^{#1}$}
\begin{thebibliography}{HGW{\etalchar{+}}17}

\bibitem[ABK98]{arbuzov1998two}
E.V. Arbuzov, A.L. Bukhgeim, and S.G. Kazantsev.
\newblock Two-dimensional tomography problems and the theory of {A}-analytic functions.
\newblock {\em Siberian Advances in Mathematics}, 8(4):1--20, 1998.

\bibitem[BGH79]{budinger1979emission}
T.F. Budinger, G.T. Gullberg, and R.H. Huesman.
\newblock Emission computed tomography.
\newblock In {\em Image reconstruction from projections. Implementaton and applications}. 1979.

\bibitem[Buk95]{bukhgeim1995inversion}
A.L. Bukhgeim.
\newblock Inversion formulas in inverse problems.
\newblock {\em Linear Operators and Ill-Posed Problems by MM Lavrentev and L.Ya. Savalev, Plenum, New York}, 1995.

\bibitem[Den23]{denisiuk2023iterative}
A.~Denisiuk.
\newblock Iterative inversion of the tensor momentum {X}-ray transform.
\newblock {\em Inverse Problems}, 39(10):105002, 2023.

\bibitem[Der23]{derevtsov2023ray}
E.Y. Derevtsov.
\newblock Ray transforms of the moments of planar tensor fields.
\newblock {\em Journal of Applied and Industrial Mathematics}, 17(3):521--534, 2023.

\bibitem[DL16]{Desai2016}
N.M. Desai and W.R.B. Lionheart.
\newblock An explicit reconstruction algorithm for the transverse ray transform of a second rank tensor field from three axis data.
\newblock {\em Inverse Problems}, 32(11):115009, 2016.

\bibitem[DPSU07]{dairbekov2007boundary}
N.S. Dairbekov, G.P. Paternain, P.~Stefanov, and G.~Uhlmann.
\newblock The boundary rigidity problem in the presence of a magnetic field.
\newblock {\em Advances in mathematics}, 216(2):535--609, 2007.

\bibitem[DS15]{derevtsov2015tomography}
E.Y. Derevtsov and I.E. Svetov.
\newblock Tomography of tensor fields in the plain.
\newblock {\em Eurasian J. Math. Comput. Appl}, 3(2):24--68, 2015.

\bibitem[DVS21]{derevtsov2021generalized}
E.Y. Derevtsov, Y.~S. Volkov, and T.~Schuster.
\newblock Generalized attenuated ray transforms and their integral angular moments.
\newblock {\em Applied Mathematics and Computation}, 409:125494, 2021.

\bibitem[Fin03]{finch2003attenuated}
D.V. Finch.
\newblock The attenuated {X}-ray transform: recent developments.
\newblock {\em Inside out: inverse problems and applications, Math. Sci. Res. Inst. Publ}, 47:47--66, 2003.

\bibitem[FOST24]{fujiwara2024inversion}
H.~Fujiwara, D.~Omogbhe, K.~Sadiq, and A.~Tamasan.
\newblock Inversion of the attenuated momenta ray transform of planar symmetric tensors.
\newblock {\em Inverse Problems}, 40(7):075004, 2024.

\bibitem[FST19]{fujiwara2019fourier}
H.~Fujiwara, K.~Sadiq, and A.~Tamasan.
\newblock A {F}ourier approach to the inverse source problem in an absorbing and anisotropic scattering medium.
\newblock {\em Inverse Problems}, 36(1):015005, 2019.

\bibitem[HGW{\etalchar{+}}17]{hendriks2017bragg}
J.N. Hendriks, A.W.T. Gregg, C.M. Wensrich, A.S. Tremsin, T.~Shinohara, M.~Meylan, E.H. Kisi, V.~Luzin, and O.~Kirsten.
\newblock Bragg-edge elastic strain tomography for in situ systems from energy-resolved neutron transmission imaging.
\newblock {\em Physical Review Materials}, 1(5):053802, 2017.

\bibitem[HS10]{holman2009weighted}
S.~Holman and P.~Stefanov.
\newblock The weighted {D}oppler transform.
\newblock {\em Inverse Problems and Imaging}, 4(1):111--130, 2010.

\bibitem[JAE{\etalchar{+}}95]{Jansson1995}
T.~Jansson, M.~Almqvist, R.~Eriksson, K.~Strahlen, G.~Sparr, H.W. Persson, and K.~Lindstrom.
\newblock Ultrasonic {D}oppler tomography-a new method to measure blood flow vector fields.
\newblock Technical report, Depanment of Electrical Measurements, Lund Institute of Technology, 1995.

\bibitem[KMM19]{krishnan2019solenoidal}
V.P. Krishnan, R.K. Mishra, and F.~Monard.
\newblock On solenoidal-injective and injective ray transforms of tensor fields on surfaces.
\newblock {\em Journal of Inverse and Ill-posed Problems}, 27(4):527--538, 2019.

\bibitem[KMS23]{kunyansky2023weighted}
L.~Kunyansky, E.~McDugald, and B.~Shearer.
\newblock Weighted {R}adon transforms of vector fields, with applications to magnetoacoustoelectric tomography.
\newblock {\em Inverse Problems}, 39(6):065014, 2023.

\bibitem[KMSS20]{krishnan2020momentum}
V.P. Krishnan, R.~Manna, S.K. Sahoo, and V.~Sharafutdinov.
\newblock Momentum ray transforms, {II}: range characterization in the {S}chwartz space.
\newblock {\em Inverse Problems}, 36(4):045009, 2020.

\bibitem[Kuc06]{kuchment2006generalized}
P.~Kuchment.
\newblock Generalized transforms of radon type and their applications.
\newblock In {\em Proceedings of Symposia in Applied Mathematics}, volume~63, page~67, 2006.

\bibitem[Lou22]{louis2022inversion}
A.K. Louis.
\newblock Inversion formulae for ray transforms in vector and tensor tomography.
\newblock {\em Inverse Problems}, 38(6):065008, 2022.

\bibitem[Mis20]{Rohit_Kumar_Mishra_moment}
R.K. Mishra.
\newblock Full reconstruction of a vector field from restricted {D}oppler and first integral moment transforms in {$\mathbb{R}^n$}.
\newblock {\em Journal of Inverse and Ill-posed Problems}, 28(2):173--184, 2020.

\bibitem[Nat01]{natterer2001mathematics}
F.~Natterer.
\newblock {\em The mathematics of computerized tomography}.
\newblock SIAM, 2001.

\bibitem[Nor89]{norton1989tomographic}
S.J. Norton.
\newblock Tomographic reconstruction of $2$-{D} vector fields: application to flow imaging.
\newblock {\em Geophysical Journal International}, 97(1):161--168, 1989.

\bibitem[Omo25]{omogbhe2025fourier}
D.~Omogbhe.
\newblock A {F}ourier approach to tomographic reconstruction of tensor fields in the plane.
\newblock {\em Journal of Mathematical Analysis and Applications}, 543(2):128928, 2025.

\bibitem[OS24]{omogbhe2024x}
D.~Omogbhe and K.~Sadiq.
\newblock On the {X}-ray transform of planar symmetric tensors.
\newblock {\em Journal of Inverse and Ill-posed Problems}, 32(3):431--452, 2024.

\bibitem[PSU13]{paternain2013tensor}
G.P. Paternain, M.~Salo, and G.~Uhlmann.
\newblock Tensor tomography on surfaces.
\newblock {\em Inventiones mathematicae}, 193:229--247, 2013.

\bibitem[RM21]{Mishra_2021}
S.K.~Sahoo R.K.~Mishra.
\newblock Injectivity and range description of integral moment transforms over m-tensor fields in $\mathbb{R}^n$.
\newblock {\em SIAM Journal on Mathematical Analysis}, 53(1):253--278, 2021.

\bibitem[Rom94]{romanov1994conditional}
V.G. Romanov.
\newblock Conditional stability estimates for the two-dimensional problem of restoring the right-hand side and absorption in the transport equation.
\newblock {\em Siberian Mathematical Journal}, 35(6):1184--1201, 1994.

\bibitem[Sha86]{sharafutdinov1986problem}
V.~Sharafutdinov.
\newblock A problem of integral geometry for generalized tensor fields on $\mathbb{R}^n$.
\newblock In {\em Doklady Akademii Nauk}, volume 286, pages 305--307. Russian Academy of Sciences, 1986.

\bibitem[Sha12]{sharafutdinov2012integral}
V.~Sharafutdinov.
\newblock {\em Integral geometry of tensor fields}, volume~1.
\newblock Walter de Gruyter, 2012.

\bibitem[SST16]{sadiq2016x}
K.~Sadiq, O.~Scherzer, and A.~Tamasan.
\newblock On the {X}-ray transform of planar symmetric 2-tensors.
\newblock {\em Journal of Mathematical Analysis and Applications}, 442(1):31--49, 2016.

\bibitem[ST15]{sadiq2015range}
K.~Sadiq and A.~Tamasan.
\newblock On the range of the attenuated radon transform in strictly convex sets.
\newblock {\em Transactions of the American Mathematical Society}, 367(8):5375--5398, 2015.

\bibitem[Ste24]{stefanov2024lorentzian}
P.~Stefanov.
\newblock The {L}orentzian scattering rigidity problem and rigidity of stationary metrics.
\newblock {\em The Journal of Geometric Analysis}, 34(9):267, 2024.

\bibitem[SU98]{stefanov1998rigidity}
P.~Stefanov and G.~Uhlmann.
\newblock Rigidity for metrics with the same lengths of geodesics.
\newblock {\em Mathematical Research Letters}, 5(1):83--96, 1998.

\bibitem[Sve12]{svetov2012reconstruction}
I.E. Svetov.
\newblock Reconstruction of the solenoidal part of a three-dimensional vector field by its ray transforms along straight lines parallel to coordinate planes.
\newblock {\em Numerical Analysis and Applications}, 5:271--283, 2012.

\bibitem[UYZ21]{uhlmann2021travel}
G.~Uhlmann, Y.~Yang, and H.~Zhou.
\newblock Travel time tomography in stationary spacetimes.
\newblock {\em The Journal of Geometric Analysis}, pages 1--24, 2021.

\end{thebibliography}
\bibliographystyle{alpha}

\end{document}